# Some Theorems on Timelike Ruled Surfaces

**Mehmet ÖNDER***
Independent Researcher, Delibekirli Village, Tepe Street, No:63, Kırıkhan, Hatay, Turkey
mehmetonder197999@gmail.com

**Abstract:** In this study, we investigate the existence theorems for timelike ruled surfaces in Minkowski 3-space $E_1^3$. We obtain a general system and give the existence theorems for a timelike ruled surface according to Gaussian curvature, distribution parameter and strictional distance. Moreover, we give some special cases such as the directirx of the surface is a geodesic, an asymptotic line, a line of curvature or a general helix.

## Timelike Regle Yüzeyler Üzerine Bazı Teoremler

**Öz:** Bu çalışmada, $E_1^3$ Minkowski 3-uzayındaki timelike regle yüzeyler için varlık teoremleri araştırılmıştır. Genel bir sistem bulunmuş ve bir timelike regle yüzeyin varlık teoremleri, Gauss eğriliği, dağılma parametresi ve striksiyon uzaklığına bağlı olarak verilmiştir. Bundan başka, yüzeyin dayanak eğrisinin geodezik, asimptotik çizgi, eğrilik çizgisi ya da bir genel helis olması gibi bazı özel durumlar verilmiştir.

## 1. Introduction

In the space, a continuously moving of a straight line generates a surface which is called ruled surface. Ruled surfaces have the most important positions and applications in the study of design problems in spatial mechanisms and physics, kinematics and computer aided design (CAD). So, these surfaces are one of the most important topics of surface theory. Because of this position of ruled surfaces, geometers have studied on these surfaces in Euclidean space and they have investigated many properties of the ruled surfaces [1-3].

Moreover, Minkowski space $E_1^3$ is more interesting than the Euclidean space. In this space, curves and surfaces have different casual Lorentzian characters such as timelike, spacelike or null (lightlike). For example, a continuously moving of a line along a curve generates a ruled surface which can be timelike, spacelike or null. Spacelike ruled surfaces are very similar to the ruled surfaces given in Euclidean 3-space $E^3$. Timelike ruled surfaces are more fascinating since there exist both timelike and spacelike curves on these surfaces. Timelike ruled surface with timelike rulings have been studied by Abdel-All, Abdel-Baky and Hamdoon [4]. Küçük has obtained some results on the developable timelike ruled surfaces in the same space [5]. Furthermore, Önder and Uğurlu have introduced Frenet frames and Frenet invariants of timelike ruled surfaces [6].

Furthermore, it is interesting to consider the existence of a timelike ruled surface in Minkowski 3-space. In global differential geometry, it is well-known that the existence and uniqueness of a surface with given first and second fundamental forms are given by Bonnet's theorem [7-9]. The theorem says that if the coefficients of these forms satisfy the Gauss equations and the Peterson-Codazzi equations, then there exists a surface, which is unique up to motions in space, for which these forms are, respectively, the first and the second fundamental forms. Of course, this is a global theorem for all surfaces and Lorentzian version of this theorem can be introduced by the similar way. But how can we give existence of a special ruled surface such as timelike ruled surface without considering its first and second fundamental forms? In this study, we try to give an answer for this question. We give some theorems for timelike ruled surfaces in Minkowski 3-space by using a similar procedure given in [10]. We obtain a general system giving the two parameter family of timelike ruled surfaces. Moreover, we give some special cases such as the directirx of the surface is a geodesic, an asymptotic line, a line of curvature or a general helix.

* Corresponding author: mehmetonder197999@gmail.com. ORCID Number of author: 0000-0002-9354-5530

## 2. Preliminaries

Let $\vec{x}=(x_1,x_2,x_3)$ and $\vec{y}=(y_1,y_2,y_3)$ be two vectors in $E^3$. The function defined by

$$<,>: E^3 \times E^3 \to IR$$
$$(\vec{x},\vec{y}) \;\to\; <\vec{x},\vec{y}> = x_1y_1 + x_2y_2 - x_3y_3$$

is called Lorentzian inner product function. The affine space $E^3$ endowed with this function is called Minkowski 3-space and denoted by $E_1^3$. In this space, an arbitrary vector $\vec{v}=(v_1,v_2,v_3)$ in $E_1^3$ can have one of three Lorentzian causal characters; it can be spacelike if $\langle \vec{v},\vec{v}\rangle > 0$ or $\vec{v}=0$, timelike if $\langle \vec{v},\vec{v}\rangle < 0$ and null (lightlike) if $\langle \vec{v},\vec{v}\rangle = 0$ and $\vec{v}\neq 0$. Similarly, an arbitrary curve $\vec{\alpha}=\vec{\alpha}(s)$ can locally be spacelike, timelike or null (lightlike), if all of its velocity vectors $\vec{\alpha}'(s)$ are spacelike, timelike or null (lightlike), respectively [11]. The norm of the vector $\vec{v}=(v_1,v_2,v_3)$ is given by $\|\vec{v}\| = \sqrt{|<\vec{v},\vec{v}>|}$.

For any vectors $\vec{x}=(x_1,x_2,x_3)$ and $\vec{y}=(y_1,y_2,y_3)$ in $E_1^3$, Lorentzian vector product of $\vec{x}$ and $\vec{y}$ is defined by

$$\vec{x}\times\vec{y} = \begin{vmatrix} e_1 & -e_2 & -e_3 \\ x_1 & x_2 & x_3 \\ y_1 & y_2 & y_3 \end{vmatrix} = (x_2y_3 - x_3y_2,\, x_1y_3 - x_3y_1,\, x_2y_1 - x_1y_2),$$

(See [12,13]).

**Definition 2.1.(See [14])** ***i) Hyperbolic angle:*** Let $\vec{x}$ and $\vec{y}$ be future pointing (or past pointing) timelike vectors in $E_1^3$. Then there is a unique real number $\theta \geq 0$ such that $<\vec{x},\vec{y}> = -\|\vec{x}\|\|\vec{y}\|\cosh\theta$. This number is called the *hyperbolic angle* between the vectors $\vec{x}$ and $\vec{y}$.

***ii) Central angle:*** Let $\vec{x}$ and $\vec{y}$ be spacelike vectors in $E_1^3$ such that $sp\{\vec{x},\vec{y}\}$ is a timelike vector subspace. Then there is a unique real number $\theta \geq 0$ such that $|<\vec{x},\vec{y}>| = \|\vec{x}\|\|\vec{y}\|\cosh\theta$. This number is called the *central angle* between the vectors $\vec{x}$ and $\vec{y}$.

***iii) Spacelike angle:*** Let $\vec{x}$ and $\vec{y}$ be spacelike vectors in $E_1^3$ that span a spacelike vector subspace. Then there is a unique real number $\theta \geq 0$ such that $<\vec{x},\vec{y}> = \|\vec{x}\|\|\vec{y}\|\cos\theta$. This number is called the *spacelike angle* between the vectors $\vec{x}$ and $\vec{y}$.

***iv) Lorentzian timelike angle:*** Let $\vec{x}$ be a spacelike vector and $\vec{y}$ be a timelike vector in $E_1^3$. Then there is a unique real number $\theta \geq 0$ such that $|<\vec{x},\vec{y}>| = \|\vec{x}\|\|\vec{y}\|\sinh\theta$. This number is called the *Lorentzian timelike angle* between the vectors $\vec{x}$ and $\vec{y}$.

## 3. Timelike Ruled Surfaces in $E_1^3$

Let $I$ be an open interval in the real line $IR$, $\vec{\gamma}=\vec{\gamma}(s)$ be a curve in $E_1^3$ defined on $I$ and $\vec{q}=\vec{q}(s)$ be a unit direction vector of an oriented timelike line in $E_1^3$. Then we have following parametrization for a timelike ruled surface $N$

$$\vec{r}(s,v)=\vec{\gamma}(s)+v\,\vec{q}(s). \tag{1}$$

In particular, if the direction of $\vec{q}$ is constant, then the ruled surface is said to be cylindrical, and non-cylindrical otherwise. The distribution parameter (or drall) of $N$ is given by

$$d=\frac{\left|\vec{\gamma}',\vec{q},\vec{q}'\right|}{\left\langle \vec{q}',\vec{q}'\right\rangle} \tag{2}$$

where $\vec{\gamma}'=\dfrac{d\vec{\gamma}}{ds}$, $\vec{q}'=\dfrac{d\vec{q}}{ds}$ (see [4,6]). If $\left|\vec{\gamma}',\vec{q},\vec{q}'\right|=0$, then normal vectors are collinear at all points of same ruling and at nonsingular points of the surface $N$, the tangent planes are identical. We then say that tangent plane contacts the surface along a ruling. Such a ruling is called a torsal ruling. If $\left|\vec{\gamma}',\vec{q},\vec{q}'\right|\neq 0$, then the tangent planes of the surface $N$ are distinct at all points of same ruling which is called nontorsal [6].

**Definition 3.1. ([6])** A timelike ruled surface whose all rulings are torsal is called a developable timelike ruled surface. The remaining timelike ruled surfaces are called skew timelike ruled surfaces. Then, from (2) it is clear that a timelike ruled surface is developable if and only if at all its points the distribution parameter $d=0$.

For the unit normal vector $\vec{m}$ of a timelike ruled surface we have

$$\vec{m}=\frac{\vec{r}_s\times\vec{r}_v}{\left\|\vec{r}_s\times\vec{r}_v\right\|}=\frac{(\vec{\gamma}'+v\vec{q}')\times\vec{q}}{\sqrt{\left\langle \vec{\gamma}',\vec{q}\right\rangle^2-\left\langle \vec{q},\vec{q}\right\rangle\left\langle \vec{\gamma}'+v\vec{q}',\ \vec{\gamma}'+v\vec{q}'\right\rangle}}.$$

Then, at the points of a nontorsal ruling $s=s_1$ we have

$$\vec{a}=\lim_{v\to\infty}\vec{m}(s_1,v)=\frac{\vec{q}'\times\vec{q}}{\left\|\vec{q}'\right\|}.$$

The plane of a skew timelike ruled surface $N$ which passes through its ruling $s_1$ and is perpendicular to the vector $\vec{a}$ is called *asymptotic plane* $\alpha$. The tangent plane $\gamma$ passing through the ruling $s_1$ which is perpendicular to the asymptotic plane $\alpha$ is called *central plane.* The point $C$ of the ruling $s_1$ where asymptotic plane is perpendicular to central plane is called *central point* of the ruling $s_1$. The set of central points of all rulings is called *striction curve* of the surface. The parametrization of the striction curve $\vec{c}=\vec{c}(s)$ on a timelike ruled surface is given by

$$\vec{c}(s) = \vec{\gamma}(s) + v_0\vec{q}(s) = \vec{\gamma} - \frac{\langle \vec{q}', \vec{\gamma}' \rangle}{\langle \vec{q}', \vec{q}' \rangle}\vec{q}$$

where $v_0 = -\dfrac{\langle \vec{q}', \vec{\gamma}' \rangle}{\langle \vec{q}', \vec{q}' \rangle}$ is called strictional distance ([6]).

**Theorem 3.1.(Chasles Theorem) ([6]):** *Let the base curve of a timelike ruled surface be its striction curve. For the angle $\mu$ between tangent plane of timelike ruled surface at the point $(s, v_0)$ of a nontorsal ruling $s$ and central plane, we have $\tan\mu = v_0 / d$ where $d$ is the distribution parameter of ruling $s$, $v_0$ is strictional distance and central point has the coordinates $(s, 0)$.*

**4. Existence Theorems for Timelike Ruled Surfaces**

Let $N$ be a timelike ruled surface in $E_1^3$ given by the parametrization

$$\vec{r}(s,v) = \vec{\gamma}(s) + v\,\vec{q}(s) \tag{3}$$

where $\vec{\gamma} = \vec{\gamma}(s)$ is directrix of $N$, $s$ is arc length of $\vec{\gamma}(s)$ and $\vec{q}(s)$ is a unit timelike vector field on $N$. The directrix $\vec{\gamma}(s)$ can be a timelike or spacelike curve. Let assume that $\vec{\gamma}(s)$ be a timelike curve. Then the Frenet formulae of $\vec{\gamma}(s)$ are given as follows

$$\begin{bmatrix} \vec{T}' \\ \vec{N}' \\ \vec{B}' \end{bmatrix} = \begin{bmatrix} 0 & k_1 & 0 \\ k_1 & 0 & k_2 \\ 0 & -k_2 & 0 \end{bmatrix} \begin{bmatrix} \vec{T} \\ \vec{N} \\ \vec{B} \end{bmatrix} \tag{4}$$

where $\vec{\gamma}'(s) = \vec{T}(s)$, $\langle \vec{T}, \vec{T} \rangle = -1$, $\langle \vec{N}, \vec{N} \rangle = \langle \vec{B}, \vec{B} \rangle = 1$, $\langle \vec{T}, \vec{N} \rangle = \langle \vec{T}, \vec{B} \rangle = \langle \vec{N}, \vec{B} \rangle = 0$, $\vec{T}, \vec{N}, \vec{B}$ are unit tangent, principal normal and binormal vectors, respectively, and $k_1$ and $k_2$ are curvature and torsion of timelike curve $\vec{\gamma}(s)$, respectively [13]. The unit normal vector $\vec{m}$ of the surface is a spacelike vector and can be given in the form

$$\vec{m} = \cos\varphi\vec{N} + \sin\varphi\vec{B} \tag{5}$$

where $\varphi = \varphi(s)$ is differentiable spacelike angle function between spacelike unit vectors $\vec{m}$ and $\vec{N}$. Let $\vec{A}$ be a spacelike unit vector in the tangent plane of the surface and be perpendicular to unit tangent $\vec{T}$. Then we can represent the timelike ruling $\vec{q}$ in the form

$$\vec{q} = \cosh\theta\vec{T} + \sinh\theta\vec{A} \tag{6}$$

where $\theta = \theta(s)$ is differentiable hyperbolic angle function between timelike unit vectors $\vec{q}$ and $\vec{T}$. It is clear that the vector $\vec{A}$ lies on the plane $Sp\{\vec{N}, \vec{B}\}$. Then we can write

$$\vec{A} = -\sin\varphi \vec{N} + \cos\varphi \vec{B} \tag{7}$$

By differentiating (5) and (7) with respect to $s$ it follows

$$\vec{m}' = k_1 \cos\varphi \vec{T} + (\varphi' + k_2)\vec{A}\,, \quad \vec{A}' = -k_1 \sin\varphi \vec{T} - (\varphi' + k_2)\vec{m} \tag{8}$$

respectively. Similarly, considering (7) and (8), the differentiation of (6) is obtained as follows

$$\begin{aligned}\vec{q}' = &\sinh\theta(\theta' - k_1 \sin\varphi)\vec{T} + \left(\cosh\theta(k_1 - \theta'\sin\varphi) - (\varphi' + k_2)\sinh\theta\cos\varphi\right)\vec{N} \\ &+ \left(\theta'\cosh\theta\cos\varphi - (\varphi' + k_2)\sinh\theta\sin\varphi\right)\vec{B}\end{aligned} \tag{9}$$

and (9) gives us

$$\begin{aligned}(\vec{q}')^2 = \langle \vec{q}', \vec{q}' \rangle = &\theta'^2 - 2k_1\theta'\sin\varphi + k_1^2(\cosh^2\theta\cos^2\varphi + \sin^2\varphi) \\ &- 2k_1(\varphi' + k_2)\sinh\theta\cosh\theta\cos\varphi + (\varphi' + k_2)^2\sinh^2\theta\end{aligned} \tag{10}$$

After these computations we can give the following theorems.

***Theorem 4.1.*** *For a timelike curve $\vec{\gamma}(s)$ as the directrix there exists a two-parameter family of timelike ruled surfaces with a given distribution parameter and a given strictional distance.*

**Proof:** The strictional distance and distribution parameter of a timelike ruled surface are given by

$$v_0 = -\frac{\langle \vec{\gamma}', \vec{q}' \rangle}{\langle \vec{q}', \vec{q}' \rangle}, \quad d = \frac{|\vec{\gamma}', \vec{q}, \vec{q}'|}{\langle \vec{q}', \vec{q}' \rangle} \tag{11}$$

respectively. Then from (6) and (9) it follows

$$v_0 = \frac{\sinh\theta(\theta' - k_1\sin\varphi)}{(\vec{q}')^2}, \quad d = \frac{\sinh\theta\left(k_1\cosh\theta\cos\varphi - (\varphi' + k_2)\sinh\theta\right)}{(\vec{q}')^2} \tag{12}$$

From (12) we have

$$\begin{aligned}\frac{d^2}{v_0^2} + 1 = &\frac{\theta'^2 - 2k_1\theta'\sin\varphi + k_1^2(\cosh^2\theta\cos^2\varphi + \sin^2\varphi)}{(\theta' - k_1\sin\varphi)^2} \\ &+ \frac{-2k_1(\varphi' + k_2)\sinh\theta\cosh\theta\cos\varphi + (\varphi' + k_2)^2\sinh^2\theta}{(\theta' - k_1\sin\varphi)^2}\end{aligned} \tag{13}$$

From (10) and (13) it follows

$$(\vec{q}')^2 = (\theta' - k_1\sin\varphi)^2\left(\frac{d^2}{v_0^2} + 1\right) \tag{14}$$

Substituting (14) in (12) we have

$$\begin{cases} \theta' = \dfrac{v_0 \sinh\theta}{d^2 + v_0^2} + k_1 \sin\varphi \\ \varphi' = -k_2 + k_1 \coth\theta\cos\varphi - \dfrac{d}{d^2 + v_0^2} \end{cases} \tag{15}$$

which are the determining equations for timelike ruled surfaces with a timelike directrix $\vec{\gamma}(s)$, given distribution parameter and given strictional distance. That proves the theorem.

***Corollary 4.1.*** *For any timelike directrix $\vec{\gamma}(s)$ as the strictional line, there exists a two-parameter family of timelike ruled surfaces with a given distribution parameter.*

**Proof:** If the directrix $\vec{\gamma}(s)$ is strictional line then $v_0 = 0$. Thus (15) becomes

$$\begin{cases} \theta' = k_1 \sin\varphi \\ \varphi' = -\dfrac{1}{d} - k_2 + k_1 \coth\theta\cos\varphi \end{cases} \tag{16}$$

that finishes the proof.

***Theorem 4.2.*** *For any timelike directrix $\vec{\gamma}(s)$ there exists a two-parameter family of timelike ruled surfaces with a given Gaussian curvature and a given angle between the tangent planes and central planes along $\vec{\gamma}(s)$.*

**Proof:** From Chasles theorem we have

$$\tan\mu = \frac{v_0}{d} \tag{17}$$

where $\mu$ is the spacelike angle between tangent plane and central plane of timelike ruled surface at the point $(s, v_0)$. Furthermore, Gaussian curvature $K$ of a timelike ruled surface is given by

$$K = \frac{d^2}{(d^2 + v_0^2)^2} \tag{18}$$

(See [15]). If we put

$$n = \sqrt{\frac{1}{K}} = \frac{d^2 + v_0^2}{d} \tag{19}$$

then $K$ defines $n$ uniquely and followings hold

$$d = n\sin^2\mu, \;\; v_0 = n\sin\mu\cos\mu \tag{20}$$

By using (20), system (15) can be given in the form

$$\begin{cases} \theta' = \dfrac{1}{n}\sinh\theta\cot\mu + k_1\sin\varphi \\ \varphi' = -\dfrac{1}{n} - k_2 + k_1\coth\theta\cos\varphi \end{cases} \tag{21}$$

that finishes the proof.

***Theorem 4.3.*** *For any timelike directrix $\vec{\gamma}(s)$ there exists a two-parameter family of general developable timelike ruled surfaces with a given strictional distance.*

**Proof:** If the timelike ruled surface $N$ is developable and not a cylinder we have $v_0 \neq 0,\ d = 0$. Then from (15) it follows

$$\begin{cases} \theta' = \dfrac{\sinh\theta}{v_0} + k_1\sin\varphi \\ \varphi' = -k_2 + k_1\coth\theta\cos\varphi \end{cases} \tag{22}$$

which shows that there exists a two-parameter family of general developable timelike ruled surfaces with a given strictional distance.

***Theorem 4.4.*** *For any timelike directrix $\vec{\gamma}(s)$ there exists a two-parameter family of timelike cylinders.*

**Proof:** If the timelike ruled surface $N$ is a cylinder then the direction of the ruling $\vec{q}$ is constant and from (9) we have

$$\begin{cases} \sinh\theta(\theta' - k_1\sin\varphi) = 0, \\ \cosh\theta(k_1 - \theta'\sin\varphi) - (\varphi' + k_2)\sinh\theta\cos\varphi = 0, \\ \theta'\cosh\theta\cos\varphi - (\varphi' + k_2)\sinh\theta\sin\varphi = 0, \end{cases}$$

that gives us

$$\theta' = k_1\sin\varphi, \quad \varphi' = -k_2 + k_1\coth\theta\cos\varphi \tag{23}$$

and (23) shows that there exists a two-parameter family of timelike cylinders.

**5. Some Special Cases**

In this section, we consider some special cases such as the directirx $\vec{\gamma}(s)$ is a geodesic, an asymptotic line or a line of curvature. Then we can give the followings.

***Theorem 5.1.*** *For any timelike directrix $\vec{\gamma}(s)$ there exists exactly one timelike ruled surface with a given Gaussian curvature on which $\vec{\gamma}(s)$ is a geodesic.*

**Proof:** Let the directrix $\vec{\gamma}(s)$ be a geodesic. Then we have $\vec{N} = \pm\vec{m}$. By considering (5), we have $\varphi = a\pi,\ (a \in Z)$. From system (21) it follows

$$\tanh\theta = \frac{nk_1}{nk_2 + 1}$$

which shows that if the directrix $\vec{\gamma}(s)$ is a geodesic, then there exists exactly one timelike ruled surface with a given Gaussian curvature.

***Theorem 5.2.*** *For any timelike directrix $\vec{\gamma}(s)$ there exists a one-parameter family of timelike ruled surfaces with a given angle between the tangent planes and corresponding central planes on which $\vec{\gamma}(s)$ is an asymptotic line.*

**Proof:** Assume that the directrix $\vec{\gamma}(s)$ is an asymptotic line. Then from (5), we have $\varphi = \pi/2$ and from the system (21) we have

$$\theta' = \frac{1}{n}\sinh\theta\cot\mu + k_1, \quad n = -\frac{1}{k_2}$$

which determines a one-parameter family of timelike ruled surfaces.

***Theorem 5.3.*** *Let the angles $\theta$ and $\mu$ be constants. Then a timelike directrix curve $\vec{\gamma}(s)$ is an asymptotic line on a timelike ruled surface $N$ if and only if $\vec{\gamma}(s)$ is a general helix.*

**Proof:** Let the directrix $\vec{\gamma}(s)$ be an asymptotic line on $N$. Then $\varphi = \pi/2$. Since the angles $\theta$ and $\mu$ are constants from (21) it follows

$$\frac{k_1}{k_2} = \sinh\theta\cot\mu$$

which is constant i.e. $\vec{\gamma}(s)$ is a general helix.

Conversely, if the angles $\theta$ and $\mu$ are constants and $\vec{\gamma}(s)$ is a general helix then from (21) we have $\varphi = \pi/2$ which shows that $\vec{\gamma}(s)$ is an asymptotic line on $N$.

***Theorem 5.4.*** *A timelike curve $\vec{\gamma}(s)$ is a line of curvature on a timelike ruled surface $N$ if and only if*

$$\varphi(s) = -\int k_2(s)ds + C \qquad (24)$$

*holds where $C$ is a constant.*

**Proof:** A curve $\vec{\gamma}(s)$ is a line of curvature if and only if surface normals generate a developable ruled surface along $\vec{\gamma}(s)$, i.e. iff

$$\left|\vec{k}', \vec{m}, \vec{m}'\right| = 0 \qquad (25)$$

Then from (8), we have $\left|\vec{T}, \vec{m}, -(\varphi' + k_2)\vec{A}\right| = 0$ which gives $\varphi' = -k_2$ and (24) is obtained.

***Theorem 5.5.*** *For any timelike curve $\vec{\gamma}(s)$ there exists a one-parameter family of timelike ruled surfaces with a given Gaussian curvature along $\vec{\gamma}(s)$ on which $\vec{\gamma}(s)$ is a line of curvature.*

**Proof:** Since $\vec{\gamma}(s)$ is a line of curvature, (24) holds and from (21) we have

$$nk_1 = \tanh\theta \sec\varphi \tag{26}$$

that finishes the proof.

***Theorem 5.6.*** *If $\theta$ is constant and timelike directrix $\vec{\gamma}(s)$ is a line of curvature on a timelike ruled surface $N$ then there exists the following relationship between the angles $\varphi$ and $\mu$*

$$\tan\varphi = -\cosh\theta \cot\mu \tag{27}$$

**Proof:** Since $\theta$ is constant and $\vec{\gamma}(s)$ is a line of curvature, from (21) and (24) we have

$$k_1 \sin\varphi = -\frac{1}{n}\sinh\theta \cot\mu, \quad k_1 \cos\varphi = \frac{1}{n}\tanh\theta \tag{28}$$

which gives (27).

## 6. Conclusions

The some theorems on the existence of timelike ruled surfaces are given by considering a timelike directrix. Furthermore, some special cases related to the directrix is introduced. Of course, one can obtain corresponding theorems for a ruled surface with null rulings or for a spacelike ruled surface.

## References

[1] Guggenheimer H. Differential Geometry. London. McGraw-Hill Book Comp. Inc. Lib. Cong. Cat. Card No. 68-12118, 1963.
[2] Karger A, Novak J. Space Kinematics and Lie Groups. Prague Czechoslovakia: STNL Publishers of Technical Lit 1978.
[3] Ravani B, Ku TS. Bertrand Offsets of ruled and developable surfaces. Comp Aided Geom Design 1991; 2(23): 145-152.
[4] Abdel-All NH, Abdel-Baky RA, Hamdoon FM. Ruled surfaces with timelike rulings. App. Math. and Comp 2004; 147: 241-253.
[5] Küçük A. On the developable timelike trajectory ruled surfaces in Lorentz 3-space $IR_1^3$. App Math and Comp 2004; 157: 483-489.
[6] Önder M, Uğurlu HH. Frenet frames and invariants of timelike ruled surfaces. Ain Shams Eng J 2013; 4: 507-513.
[7] Blaschke W, Leichtweiss K. Elementare Differentialgeometrie. Springer, 1973.
[8] Do Carmo M. Differential geometry of curves and surfaces. New Jersey. Prentice-Hall, 1976.
[9] Struik DJ. Lectures on classical differential geometry. New York. Dover Publications, Inc. 1988.
[10] Kamenarovic I. Existence Theorems for Ruled Surfaces in the Galilean Space $G_3$. Rad Hrvatske Akad. Znan. Umj. Mat 1991; 10(456): 183-186.
[11] O'Neill B. Semi-Riemannian Geometry with Applications to Relativity. London. Academic Press, 1983.
[12] Uğurlu HH, Çalışkan A. Darboux Ani Dönme Vektörleri ile Spacelike ve Timelike Yüzeyler Geometrisi. Manisa. Celal Bayar Üniversitesi Yayınları. Yayın No: 0006, 2012.
[13] Walrave J. Curves and Surfaces in Minkowski space. PhD. thesis, K.U. Leuven: Fac. of Science, Leuven, 1995.
[14] Ratcliffe JG. Foundations of Hyperbolic Manifolds. Springer, 2006.
[15] Kazaz M, Özdemir A, Güroğlu T. On the determination of a developable timelike ruled surface. SDÜ Fen-Edebiyat Fakültesi Fen Dergisi (E-Dergi) 2008; 3(1): 72-79.